\documentclass[twoside]{article}
\usepackage[section]{algorithm}
\usepackage{amsthm,amsmath,amssymb,amsmath,amsfonts,graphicx,fancyhdr,algorithmic,longtable,enumitem,epsfig}
\usepackage{epstopdf}
\usepackage{bibcheck}
\usepackage{listings}
\usepackage{hyperref}
\usepackage{appendix}

\usepackage{multicol}
\def\enddemo{\qed \endtrivlist}
\expandafter\let\csname enddemo*\endcsname=\enddemo

\def\qedsymbol{\ifmmode\bgroup\else$\bgroup\aftergroup$\fi
  \vcenter{\hrule\hbox{\vrule height.6em\kern.6em\vrule}\hrule}\egroup}
\def\qed{\ifmmode\else\unskip\nobreak\fi\quad\qedsymbol}

\newtheorem{algo}{Algorithm}[section]

\newtheorem{fig}{Figure}[section]

\font\sst=cmtt8

\font\ssi=cmti8

\font\srs=cmr6

\newcommand{\eqdef}{\overset{\mathrm{def}}{=}}

\title{\bf Vilin: Unconstrained Numerical Optimization Application}

\author{\frenchspacing
\bf Marko Miladinovi\' c$^{1}$\footnote{Corresponding author}, Predrag \v Zivadinovi\' c$^{2}$, \\
$^{\srs 1,2}${\ssi University of Ni\v{s}, Faculty of Sciences and Mathematics, Department of Computer Science}\\
{\ssi Vi\v segradska 33, 18000 Ni\v s, Serbia.
\footnote{The authors gratefully acknowledge support from the research project 144011 of the Serbian Ministry of Science.}
} \\
{\ssi {\ssi E-mail}:} $^{\srs 1}${\sst markomiladinovic@gmail.com}\ $^{\srs 2}${\sst predrag.zivadinovic@pmf.edu.rs}}
\date{} 
\begin{document}

\maketitle

\begin{abstract}
We introduce an application for executing and testing different unconstrained optimization algorithms.
The application contains a library of various test functions with pre-defined starting points.
A several known classes of methods as well as different classes of line search procedures are covered.
Each method can be tested on various test function with a chosen number of parameters. Solvers come
with optimal pre-defined parameter values which simplifies the usage. Additionally, user friendly
interface gives an opportunity for advanced users to use their expertise and
also easily fine-tune a large number of hyper parameters for obtaining even more optimal solution.
\smallskip

\indent This application can be used as a tool for developing new optimization algorithms (by using simple API),
as well as for testing and comparing existing ones, by using given standard library of test functions.
Special care has been given in order to achieve good numerical stability of all vital parts of the application.
The application is implemented in programming language Matlab with very helpful gui support.

\frenchspacing \itemsep=-1pt
\begin{description}
\item[] AMS Subj. Class.: {68N01, 90C30, 90C06}
\item[] Keywords: Software, Line search, Gradient descent methods, Newton method, Quasi Newton methods,
Modified Newton methods, Conjugate gradient methods, Trust region methods.
\end{description}
\end{abstract}

\section{Introduction}\setcounter{equation}{0}
The theory of unconstrained optimization is widely used and is, also, active field of research in the different areas such as
military, economy, finance, mathematics, and computer science. Many papers has been published in this field in the recent years
\cite{Andrei1}, \cite{Andrei2}, \cite{Jamil}, \cite{Jian}, \cite{Kafakia}, \cite{Rao}, \cite{Wang}.
\smallskip

During the long history of unconstrained optimization the researchers have made a lot of efforts to create and invent a large variety of different test functions.
It is clear that the well  designed  test  problems  are  very  helpful  in clarifying the new methods ideas and mechanisms. Also, a reasonably large set of
test problems need to be presented in order to get a clear conclusion about the hypothesis used in proving the  quality  of  the  algorithm.
These conclusions are made based on comparison of local and global convergence, complexity of the algorithm, CPU time, number of iterations and other
important features obtained at an experimental level. Additionally, numerical stability and rate of convergence as a very important features can be
used in making decisions and conclusions.
\smallskip

Bearing all this in mind there are some questions that may arise and concern many researchers in this field.
Very similar questions are already presented in \cite{Bongartz}
\begin{enumerate}[parsep=-1pt]
\item[$\bullet$] I want to test and compare my solver with other existing solvers on standard test set. Is this
test set available?
\item[$\bullet$] I know which problem(s) I want to solve. Which solvers are available and how can I use them?
\item[$\bullet$] The solver that I want to use is not available. Is there an easy way of building a suitable 
interface to test this solver on some existing or other specific class of problems?
\item[$\bullet$] Is there any stable and reliable environment(application) with standard test set and standard solvers?
\item[$\bullet$] If application, with standard library of tests functions, exists is there an easy way to include the new test function?
\item[$\bullet$] If standard application (environment) exists is there an easy way to include the new or modification of
existing solver?
\end{enumerate}

In order to give an answer on some of this questions, researches developed different test functions.
An important source of goal functions is the CUTE collection established and published by Bongartz
et al in \cite{Bongartz}.
Later, as an extension and improved version of this collection, Gould et al published another collection
and name it as CUTEr \cite{cuteR}.
Some other test functions are presented by More at al in \cite{More} and Himmelblau in \cite{Himmelblau}.
The most important from this collection, as well as  some other from different sources, are well studied,
collected, systematized and presented in algebraic form
by Andrei in \cite{AndreiTest}. It is important to say that a lot of researchers have contributed to the preparation of this very
nice collection, as already mentioned. This collection belongs to the group of {\em artificial unconstrained optimization test problems}.
Each test function is accompanied with appropriate starting point $x_0$ which is very important input for testing and comparing
performances of different methods.  By collecting and standardizing the test functions in to the library such as the one given
in \cite{AndreiTest} the process of testing and comparing different methods of unconstrained optimization is significantly facilitated.
As we already emphasized there is a clear need of standardizing and simplifying the process of comparing different and similar methods.
It is inevitable that, during the process of researching and developing some new methods
or modifications of existing methods, the researchers concern themselves with problem
of testing and comparing with other methods.
By having a nice test functions collection, such as one given by \cite{AndreiTest}, one part of this process is completed.
\smallskip

Another very important part is to create (develop) a framework (application) as a very useful tool for testing a number of different
solvers onto the various number of standardized test functions.
One of the first and most important package which deals with this problem is very well know LANCELOT package
introduced by Conn at al in \cite{Lancelot}. Namely, this package present the library for large scale nonlinear optimization
problems written in programming language Fortran. The authors claimed that big contribution has been made in the optimization community, concerning
both theoretical aspects as well as software for constrained and unconstrained large scale problems. Additionally, numerical
experiments and analysis of the methods from LANCELOT package have been presented in \cite{LancelotExp}. Also, the algorithmic options
of the methods from the library have been discussed.
Next platform that is very powerful is well known GALAHAD library which contains $90$ packages for
solving large-scale nonlinear optimization problems. This platform introduced by Gould at al \cite{Gould} is written in programming language Fortran.
Package also contains the updated versions of the nonlinear programming package LANCELOT.
Over the years these platforms have been heavily used, especially  in the optimization community, for solving variety of large scale problems.

\smallskip
Besides these libraries there is a clear need for the appearance of some powerful
and user friendly platform with user friendly graphical interface. This graphical interface can further simplify
and accelerate the whole process of testing and comparing. The user friendly platform can give the answers on the remainder of
the mentioned questions that concerns researchers from this field. Besides the fact that this application should contain the standard solvers
as well as the standard test functions, it should be very easy to operate with and should possess a property of easily adding new solvers and test
functions. After providing mentioned functionality some of the powerful benchmark for testing and comparing different methods can be applied.
One of the most powerful and most widely used metric for benchmarking and comparing optimization software is performance profile metric
introduced by Dolan and More in \cite{Dolan}.

\smallskip
In this paper we present a Matlab \rm{gui} application which contains standard unconstrained solvers and standard test functions taken from
\cite{AndreiTest}. We believe that most of the concerns that we emphasise and the possible issues that arise are
covered and resolved by this user friendly platform.

\smallskip

\textit{Overview}: The remainder of the paper is organized as follows. In section \ref{ArchDesign} architecture and design of the application are shown,
its API and explanation of the concept. This section is separated in three subsections. The first one presents the concept of test
functions and the procedure of adding new functions. 
The given groups of unconstrained methods as well as the basic explanation of the code organization
are covered in the second subsection. In the third subsection directions on how different line search methods
should be used (and expanded) are given.  Section \ref{ApplicationOverview} presents the detailed application overview,
in particular, its graphical frontend and also nice explanations of how to use it. Additionally, as a separate subsection,
some important notes and instructions for application users are provided.
Conclusions and possible further steps related to the upgrade of the application are presented in section \ref{Conclusions}.
Finally, we provide another section under appendix \ref{MethodsList}, which gives the full list of available methods and line 
searches with some very brief explanations as well as reference to the original papers. 

\section{Architecture and design}\label{ArchDesign}\setcounter{equation}{0}

In this section we present the platform for testing and comparing different standard unconstrained solvers on some standard
test functions taken from \cite{AndreiTest}. By this application with nice graphical interface we try to simplify the whole
process of testing solvers and comparing different methods on standard test set.
\smallskip

Complete application is written in Matlab programming language. Motivation for this decision is in simple and powerful notation of the language,
very good plotting support and decent GUI capability. One of the goals was to provide framework for easy and straightforward developing and
testing of new methods and functions. Matlab is good choice because it is ubiquitous and well known in optimization as well as numerical
community. Three parts of the application are of the most interest, objective functions, optimization methods and line search methods.
In the rest of this section main properties of each part will be covered as well as the process of extending them.
\smallskip

We would like to point out that the during of the process of implementing the methods, line searches and objective functions
the special attention was put on code efficiency and numerical stability. Namely, each method is written according to the original
paper as well as some possible improvements which appeared later.
The field of improving the methods is very dynamic and active and lot of researches are still working on code stability and
efficiency. Therefore, some of the given method implementations may be slightly worse than the state of the art in numerical community.
We strongly believe that, in future releases, we will have even better and more efficient method implementations.
\smallskip

We provided the source code for the application which is publicly available and can be downloaded on the following link
\href{https://github.com/markomil/vilin-numerical-optimization}{\underline{vilin}}.

\subsection{Objective functions}\label{ObjectiveFunctions}

In optimization, the main goal is to find extremum of appropriate function. This function is referred to as objective function. A lot of effort has been
made to collect functions for testing different optimization methods. Most of them (obtained from various sources) are part of the presented application.
Mostly, functions are taken from unconstrained test collection \cite{AndreiTest}.
This collection is very important in the field of unconstrained optimization. As a confirmation of this statement, in recent
years a large number of articles which cited the paper \cite{AndreiTest}, have been published, see for example
\cite{Ahookhosh}, \cite{Andrei1}, \cite{Jamil}, \cite{Jian}, \cite{Kafakia}, \cite{Wang}. The presented test functions are well prepared and collected.
For each test function it's algebraic representation is given. Functions are presented in extended or generalized form.
The main difference between these forms is that while the problems in extended
form have the Hessian matrix as a block diagonal matrix, the generalized forms have the Hessian as a multi-diagonal matrix (tridiagonal, pentadiagonal etc). Generally,
the problems given in generalized forms are slightly more difficult to solve (require more iterations and CPU time).
\smallskip

During the process of minimization each optimization method requires function value,
numerical gradient and/or numerical Hessian at some point to be calculated. Thus, test function need to be presented by appropriate matlab function which return
mentioned numerical values (objects). The set of three objects or some subset of these important numerical values at
given point $\bold x_k$ need to be provided. Below is shown prototype of a goal function from application.

\begin{small}
\begin{verbatim}
function [ outVal, outGr, outHes ] = ExtRosenbrock(x_k, VGH)
\end{verbatim}
\end{small}

Each objective function is implemented following the same pattern which is standard in numerical
optimization community.
Function accepts two arguments and returns three already mentioned objects values.
First input argument is point $\bold x_k$ at which one wants to evaluate the mentioned objects values,
while second argument is vector ${\texttt{VGH}}$ of length $3$, which represent the flag for computing function value, gradient and Hessian.
This input argument is formed
in the following way: if function value is needed then $\texttt{VGH}(1) = 1$; if gradient is required then $\texttt{VGH}(2) = 1$; if Hessian need to be
computed then $\texttt{VGH}(3) = 1$. Otherwise these values are set to zero. Each function returns three parameters
computed value, gradient and Hessian respectively.
By default, if any output object is not implemented (for chosen goal function) and need to be used by some optimization method
an error will be thrown and helpful message will be displayed in user interface so problems can easily be resolved.
\smallskip

In order to have a unified template for testing and comparing different methods each objective function is accompanied with appropriately chosen starting point.
Starting points and given objectives for this purpose are taken from \cite{AndreiTest}.
Default starting points for each test functions are predefined. Their implementations are stored in file '\texttt{Util/StartingPointGenerator.m}'.
Starting points definitions are located at the beginning of the file.
\newpage

\vskip 10pt
{\bf \noindent Adding new test function}
\vskip 10pt

In our application the test functions are located in folder '\texttt{Functions/MultiDimensional}'.
For easier adding of new functions template with appropriate code is provided at '\texttt{Functions}' directory, see
'\texttt{NewFunctionTemplate.m}'. In order to add new test
function one have to change the given template and to save it in '\texttt{Functions/MultiDimensional}' directory.
After running the application, added function will become visible in function popup.
\smallskip

In order to provide the implementation of the new starting point (accompanied to new added test function) new entry at the end of
the list should be added in file '\texttt{Util/StartingPointGenerator.m}'.
Required dimension is passed to the file, also. 
Starting points are computed by following the specific rule, so logic for generating must be implemented.
However, in most cases this is already done and usually one just need to call appropriate function (which represent the starting point generator)
and to pass dimension
of the starting point. Implemented rules are located in the second part of the same utility file. If needed the code for new rules can be easily provided.

\subsection{Methods explanation}\label{Methods}

This application is built as a framework for testing, comparing and developing optimization methods having in mind simplicity
and productivity. Thus, special effort has been put in order
to achieve these goals. As a starting reference for work many basic and well known methods have been implemented and are now available inside the application.
However, this list is incomplete and it should be expanded in the future. Presented methods are iterative and most of them follows the general scheme
\begin{equation}\label{MethodsGeneralForm}
\bold x_{k+1} = \bold x_k + t_k \bold d_k, \quad k = 0, 1,\ldots
\end{equation}
Point $\bold x_k$ denotes approximation of function minimum at $k$-th iteration, $\bold d_k$ denotes search
direction while $t_k$ represents step size in direction $\bold d_k$.
\smallskip

For the purpose of easy managing and using, methods are grouped by their nature and similarity in the following six groups:
\begin{itemize}[parsep=-4pt]
  \item[$\bullet$]Gradient Descent
  \item[$\bullet$]Newton
  \item[$\bullet$]Conjugate Gradient
  \item[$\bullet$]Modified Newton
  \item[$\bullet$]Quasi Newton
  \item[$\bullet$]Trust Region
\end{itemize}

Some short explanations about the properties of each of the methods group as well as the list of available and given
methods are presented in the Section \ref{MethodsList}.
Below is shown prototype of a method from application.

\begin{small}
\begin{verbatim}
function [ fmin, xmin, iterNum, cpuTime, evalNumbers, valuesPerIter ] =
          DaiYuan(functionName, methodParams)
\end{verbatim}
\end{small}

Each method has two input arguments, name of the objective function and method parameters. Clearly, function name points to
objective function, while '\texttt{methodParams}'
holds parameters for method. These parameters are encapsulated in class '\texttt{MethodParams}', located
in folder '\texttt{Utils}', and passed to method before execution.
Parameters are initialized from values that are filled in through application interface.
Also, name of chosen line search and its parameters values are also passed through
this argument and segregated inside method's body before passing to specific line search method, if any is used.
\smallskip

Methods have six output arguments. First is minimal value of the objective function that is found by the method. Second one is a point at which
function minimum is reached. These two arguments are obtained as a result of applying method onto the objective function.
Other arguments are used for displaying important properties and later examining performances of methods.
Third returned argument represents required number of iterations, fourth is total CPU time elapsed.
The following returned argument holds evaluation numbers of function value, gradient and Hessian, respectively.
This is also encapsulated in separate class called '\texttt{EvaluationNumbers}' and
located in folder '\texttt{Util}'. Last argument is used for plotting function values and gradient norm
values through iterations. This argument holds all necessary values in class '\texttt{PerIteration}'.

\subsubsection{Adding new method}

Methods are located in folder '\texttt{Methods/MultiDimensional}' and each of mentioned groups has its separate subfolder.
Methods are written in similar manner as functions, following appropriate template which can be found in the same folder as methods,
see '\texttt{NewMethodTemplate.m}'. Therefore, workflow for adding new methods is easy as the procedure for adding new functions.
One only needs to change template located in folder '\texttt{Methods/MultiDimensional}' and to put it in appropriate subdirectory.

\subsection{Line search methods explanation}\label{LineSearchMethods}

As already mentioned most of the methods follow the general iterative form \eqref{MethodsGeneralForm}.
Clearly, two problems that arise are how to find search direction $\bold d_k$ and how to compute appropriate step size $t_k$.
In order to determine search direction $\bold d_k$, different methods are provided, see Section \ref{MethodsList}.
Ones search direction $\bold d_k$ is computed the problem of finding step size 
can be reformulated as finding $t_k$ such that
\begin{equation} \label{LineSearchCondition}
\Phi(t_k) < \Phi(0), \quad {\rm where} \quad \Phi(t) \eqdef f(\bold x_k + t \bold d_k) = f(\bold x_{k+1})
\end{equation}
is one dimensional function. In order to determine the step size parameter different line search procedures can be applied.
There exists several different algorithms that satisfies number of distinct line search rules.
In application we covered the most important implementations of different line search rules. The list of available line searches, some
brief explanation of each line search and its parameters are presented in the section \ref{lineSearchList}.
Below is shown prototype of a line search method.

\begin{small}
\begin{verbatim}
function [ outT, outX, outVal, outGr, evalNumbers ] = ApproxWolfe( functionName, params)
\end{verbatim}
\end{small}

From above definition it can be seen that, similar as main methods, line search methods take function name and parameters '\texttt{params}' as inputs.
Params holds some internal parameters values which can be set from the application interface or inherit default predefined values. They  are encapsulated and explained
in file '\texttt{LineSearchParams.m}' located in folder '\texttt{Util}'. Object of this class is constructed in main method's body, then passed to line search method.
Line searches return five output values, desired step size $outT$, new point $outX$ (after computing step size), current function value ($outVal$),
current gradient value ($outGr$), as well as evaluation numbers (explained earlier) that are included in main method's statistics.

\subsubsection{Adding new line search method}
The line search methods are located in folder '\texttt{Methods/MultiDimensional/LineSearch}'.
Adding new line search methods is supposed to be as easy as adding new functions or methods.
Therefore, one only needs to change template called '\texttt{NewLineSearchTemplate.m}' located in same folder as
other methods and line searches and to put it to appropriate subdirectory '\texttt{MultiDimensional/LineSearch}'.

\section{Application overview}\label{ApplicationOverview}\setcounter{equation}{0}

In this section we present the application's overview. Namely, we go in to the details of one of the main part of this platform which is a
very nice and useful graphical interface. Also, some explanations of the application functionality and capability
from the user point of view are provided. Additionally, some very important notes about the efficient application usage are presented.

\subsection{GUI overview}

In order to start using this platform one only need to run a main function '\texttt{vilin.m}' in root folder. After running the code the main application
window will popup, see figure \ref{numOptApp}. The application window consists of two main parts.
On the left half of the screen one can find controls that accept users input
and the main button of the application. Two main figures cover the most of the right hand part of the screen.
The results of running particular method on chosen objective function are shown in the lower part of the right side under the panel
named '\texttt{Results}', see figure \ref{numOptApp}.

\begin{center}\small
$\includegraphics* [scale=0.35]{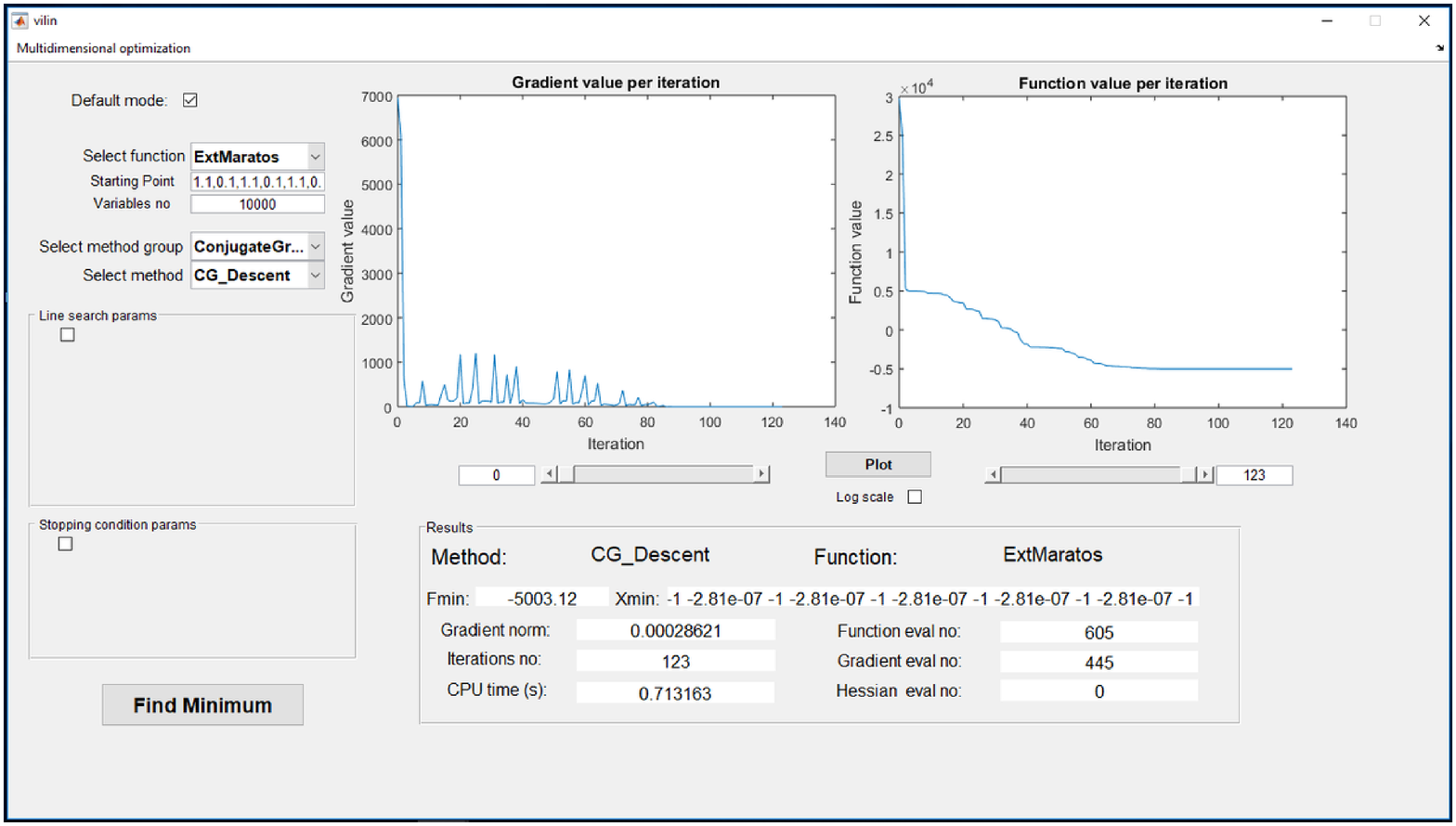}$
\begin{fig}\label{numOptApp}
Numerical optimization application.
\end{fig}
\end{center}

\vskip 5pt
{\bf \noindent Input controls}
\vskip 5pt

The left part of the application window consists of several popup boxes, several edit boxes,
two panels and the main application button '\texttt{Find Minimum}'.
\smallskip

First popup box is used for choosing the test function, the objective function on which one
wants to run a specific optimization method. One of many different test functions can be chosen.
Below the function popup box two edit boxes are located. They are related to starting point $x_0$ for chosen objective function
as well as it dimension refers as '\texttt{Variables no}'. By default, for each test function starting point value
is predefined, see \cite{AndreiTest}. It is important to note that each starting point follows some rule.
Thus, the changes of the dimension produces the changes of the starting point values, immediately, according to the given rule.
Detailed instructions on adding new functions and therefore adding starting points and possible new rules are covered
in section \ref{ObjectiveFunctions}.
\smallskip

Second popup box serves for choosing an appropriate method group. Namely, all provided solvers are divided in six groups: Conjugate Gradient,
Gradient Descent, Modified Newton, Newton, Quasi Newton and Trust Region. Each group consists of one or more methods. Therefore, in order to choose a
method one first has to select the group that its belongs. Later, appropriate method from the chosen group can be selected by
the third popup box.
\smallskip

First panel is '\texttt{Line search params}' panel whose objects  are hidden by default. This panel is visible only for those solvers that use
line searches to compute step-size, otherwise is invisible. In order to enable parameters for tuning one just need to mark the given checkbox.
Once the checkbox is marked several edit boxes as well as one popup box will appear. All these objects refers to the procedure for
choosing line search rule as well as tuning appropriate parameters values. Line search can be chosen by popup box while
appropriate parameters can be tuned by using edit boxes. However, most of the solvers come with default and optimal set-up which
includes pre defined line search procedure and optimal hyper parameters values. Additionally, the advanced user can use knowledge and expertise to try to
find even more optimal settings for the specific problem by changing some of these input objects values and possible obtain more optimal results.
Parameters that belongs to this panel are the following parameters $\beta$, $\rho$, $\sigma$, $t\_Init$ and $M$.
The detailed explanation about the parameters  as well as the
parameters boundaries are covered in section \ref{lineSearchList}.
\smallskip

Second panel named '\texttt{Stopping condition params}' covers the hyper parameters for controlling the stopping criteria of the
iterative process for the chosen solver. Similarly as for previous panel the objects of this panel are invisible
by default. Namely, the panel contains three edit boxes that cover the maximal number of iterations,
algorithm accuracy (achieved by the gradient norm) and working precision. To summarize, each of the given methods is iterative,
thus it is possible to choose maximal number of iterations as well as method precision.
Execution will stop when method reaches maximal number of iterations or when gradient norm becomes lower
then '\texttt{Epsilon}' value. Additionally, for the sake of numerical stability '\texttt{Work precision}' is added to ensure that
function value changes significantly between adjacent iterations.
Therefore, the termination criteria is as follows:
\begin{equation}
{\rm numIter} \geq {\rm maxIterNum}, \quad {\rm or} \quad \| \bold \nabla f(\bold x_k) \| \leq
\varepsilon, \quad  {\rm or} \quad \frac{|f(\bold x_{k-1})-f(\bold x_k)|}{1+|f(\bold x_k)|} \leq {\rm workPrec}.
\end{equation}

Besides the explained objects, there is one additional checkbox that need to be explained, named '\texttt{Default mode}'. 
Namely, this checkbox is marked by default which means that for each method or method group appropriate set up will be
chosen which best fit the selected method. For example, if '\texttt{CG\_Descent}' conjugate gradient method is chosen
the appropriate '\texttt{ApproxWolfe}' line search rule will be preselected which comes with appropriate pre-defined parameter
values. This means that under active '\texttt{Default mode}' appropriate line search will be preselected (if necessary)
as well as some predefined hyper parameters values. All this values and choices are taken from the authors original
papers or from some later articles that improves the original methods.
The main purpose of providing default mode option is the ability for each method to
have predefined optimal settings that goes with the method which simplifies the application usage.
Nevertheless, the advanced users will have the opportunity to do some additional experiments
(with different line searches and different parameters values) by unchecking the '\texttt{Default mode}'.
Therefore, if '\texttt{Default mode}' checkbox is not checked changing the selection of the method group or method
won't affect in selection of the line search method.

\vskip 10pt
{\bf \noindent Output objects}
\vskip 10pt

Once the process of choosing appropriate test function as well as the solver is finished the process of function minimization starts by clicking
on button '{\rm Find Minimum}', see figure \ref{numOptApp}. When computation is
finished results become visible on the right hand side of the screen which is intended for displaying the output results.
\smallskip

The upper part of output area (right hand side of the screen) is occupy by two figures. The first
one shows how norm of the function gradient changes through iterations, while the changes of the function value through iterations
are presented in the second one. Given information are of great importance for analyzing methods complexity and their
behaviour on different test functions. Sliders below figures can be used to show smaller parts of plots,
between two arbitrary iterations. This is very useful when computation takes large number of iterations and when only small portions of
them have to be looked at and analyzed precisely. These small portions can give, additional, very useful information about the convergence of the method
and significantly improve one's knowledge of the nature of the solver. Additionally, checkbox for displaying results in logarithmic scale is provided.
By marking this checkbox, instead of original, figure displays logarithm of original values which gives more useful information about the iterative process
of the solver.
\smallskip

In the lower part of the output area the output results, after applying the particular solver on the goal function, are shown.
These output values are gathered and displayed on panel named '\texttt{Results}', see figure \ref{numOptApp}.
The obtained function minimum '\texttt{Fmin}' and the point '\texttt{Xmin}', in which the minimum is reached, are given in first row.
Additionally, in the left column, gradient norm achieved in the last iteration, total number of iterations and total time spent in computing
'\texttt{CPU time}' are presented. In the right column one can see how many times
numerical values of function, gradient and hessian matrix are computed . These numbers represent standards on measuring computational complexity of the
selected method. 
\smallskip

To conclude, the total performance of the method is presented by the given output data. These data matter in comparing different methods and
give enough information in analyzing and finding weak spots in process of developing new or improving existing methods.

\subsection{Important notes and brief instructions}\label{notes}
Here we present some short notes and instructions about the way application should be used.

\begin{itemize}[parsep=-1pt]
\item As we already mentioned, for the purpose of easy usage and user friendly environment, the application
provides so called '\texttt{default mode}' settings. This mode implies optimal parameters set-up for each solver and thus
fast and easy way of optimizing goal functions.
\item In addition, by unchecking the '\texttt{default mode}' checkbox, the advanced users can manually do some
fine tuning of the input parameters values.
Namely, even though each solver comes with predefined and optimal set-up, application gives an opportunity for combining
different line searches as well as tuning other input parameters values for the chosen solver. This is very useful in
the process of creating new optimization methods as well as improving the existing ones.
It is not recommended to try some arbitrary selection that is not supported by the optimization theory.

\item The users should be very careful with combining line search procedures with chosen solver as well as with additional parameters
tuning. Namely, for most line search procedures there exist explicit parameter boundaries that should be satisfied.
Also, in order to converge some solvers need to be accompanied by appropriate line search algorithms.
The user should follow the general rules known in optimization community and that are also covered in original
papers or some monograph, see for example \cite{FletcherBook, Nocedal}.

\item To conclude, for those users with no theoretical background the '\texttt{Default mode}' is the best choice.
Otherwise, the application provides the ability of easy testing of different combinations of methods and line searches
as well as different parameter values tuning.
\end{itemize}


\section{Conclusions and Future Work}\label{Conclusions}\setcounter{equation}{0}
An application for benchmarking and developing optimization methods has been introduced.
The main idea was to create a tool that will be used in academia for developing new methods
as well as for presenting and teaching well known solvers to the students in the field of unconstrained numerical optimization.
Additionally, people working in industry may find this application useful for testing and comparing different methods and
deciding which one is the most applicable for their practical problems.
For achieving these goals some common requirements need to be satisfied.
One of the most important requirements of every optimization software is numerical stability, so
special care has been taken in order to provide good code efficiency and numerical stability of the solvers.
Other requirements, that may arise and concern researchers in this field, are given in the form of questions and are already presented
in the introduction. We believe that the goal is achieved by building this application
with very helpful graphical user interface while exposing powerful, yet simple, API for extending its capabilities.
\smallskip

\noindent The main advantages of the introduced platform are
\begin{enumerate}[parsep=-1pt]
\item[$\bullet$] Standard test functions library availability
\item[$\bullet$] Simple procedure for adding new test problems
\item[$\bullet$] Standard solvers library availability
\item[$\bullet$] Simple and standardize way for adding new and modifying existing solvers
\item[$\bullet$] Simplify procedure for comparing different solvers with respect to the many important features
\item[$\bullet$] Very nice graphical support that gives more information about the nature of the solver with respect to the chosen test function
\item[$\bullet$] Simple and fast ability of combining different line searches with appropriate methods
\item[$\bullet$] Simple and fast procedure of experimenting with different hyper parameters values
\end{enumerate}

To conclude, this platform is unique and so far, as we know, there is not a similar application available.
We strongly believe that this platform can do much to simplify the researcher's work and to provide them with
new capabilities that until now weren't available.
\smallskip

In the future we plan to add support for creating various test scenarios. This includes possible selection of some subset of test functions. This selection can be useful as a
comprehensive benchmark of method properties, as well as for generating reports on such tests. With this in mind web based graphical
interface arise as another idea that we have in plan for improving our application capabilities. This web based interface can be very handy because users
will be allowed to deploy application on server and to remotely access it and run different test examples. Additionally, with growing interest in research and
applications of deep learning, whose core part is in solving optimization problems, an idea for adding support for testing methods
on custom neural network architectures is considered and will be inspected in detail in near future.
\smallskip

The complete source code is available and can be downloaded on the following
\href{https://github.com/markomil/vilin-numerical-optimization}{link}.

\begin{appendices}
\section{Appendix: Methods overview}\label{MethodsList}\setcounter{equation}{0}

In this section we present the list of available methods which belongs to the specific method group. Some very shorts explanations are provided.
For detailed information about the specific method see the original paper or some relevant book or monograph, such as \cite{FletcherBook},
\cite{Nocedal} and \cite{SY}.
For the sake of simplicity we introduce the following notations:
$$\bold g(x) = \nabla f(\bold x), \quad \bold g_k = \nabla f(\bold x_k), \quad \bold s_k = \bold x_{k+1} - \bold x_k,
    \quad \bold y_k = \bold g_{k+1} - \bold g_k, \quad G_k = \nabla^2 f(\bold x_k).$$ 

\subsection{Gradient descent}

This group of methods represent the first-order iterative optimization algorithm group.
Original method takes steps proportional to the negative of the gradient of the function at the current point,
also known as gradient descent method. Additionally, group also covers some modifications of original gradient descent.
Each of the method follows simple iterative rule
\begin{equation}\label{gradient_descent}
\bold x_{k+1} = \bold x_k - t_k \bold g_k, \quad k = 0,1,\ldots
\end{equation}
or some slight modification. Three different methods from this group are covered.

\vskip 8pt
{\bf \noindent Gradient descent with line search}
\vskip 5pt
Gradient descent also known as steepest descent algorithm is one of the simplest and most famous methods in the theory of
unconstrained optimization. It is introduced by Cauchy as a method for solving the system of linear equation \cite{Cauchy}.
Gradient descent (given by \eqref{gradient_descent}) has linear rate of convergence, but it's convergence
is inferior to many other methods.
For poorly conditioned convex problems, gradient descent increasingly 'zigzags' as the gradients point nearly orthogonally to
the shortest direction to a optimal point which significantly decrease the convergence. This method can be combined with various
line search algorithms, covered in section \ref{lineSearchList}.
\vskip 8pt
{\bf \noindent Barzilai-Borwein}
\vskip 5pt
Barzilai-Borwein is two-point step size gradient method, originally introduced by its authors Barzilai and Borwein \cite{Barzilai}.
Main idea was to determine two-point step sizes for the steepest descent method by approximating the secant equation as follows
\begin{equation}\label{BBgama}
t_{k+1}^{BB}=\frac{\bold s_k^T \bold y_k}{\bold y_k^T \bold y_k}, \quad k = 0,1,\ldots
\end{equation}
Authors claim that the algorithm achieves better performance and cheaper computation than the classical steepest descent method.
Later, in order to make this method globally convergent and thus more efficient, Raydan in his paper \cite{Raydan} accompanied a nonmonotone
line-search procedure proposed by Grippo at al \cite{Grippo}. The resulting method was imposed as an efficient solver
for large scale unconstrained minimization problems.
\vskip 8pt
{\bf \noindent Scalar Correction}
\vskip 5pt
Another two-point step size method, named Scalar Correction method, is introduced by Miladinovi\'c at al \cite{Miladinovic}.
The initial trial step-length is determined from the secant equation as well as Hessian inverse approximation by an appropriate scalar matrix
\begin{equation}\label{scalar_correction}
t_{k+1}^{SC} = \left\{
\begin{array}{ll}
    \frac{\bold s_k^T \bold r_k}{\bold y_k^T \bold r_k}, & \bold y_k^T \bold r_k > 0 \\
    \frac {\|\bold s_k\|}{\|\bold y_k\|}, & \bold y_k^T \bold r_k \leq 0.
\end{array}
\right.
\end{equation}
where $\bold r_k = \bold s_k- t_k\bold y_k$. In order to get globally convergent algorithm, nonmonotone line search procedure is
accompanied. The reported numerical results indicate improvements of the performance with respect to the global Barzilai-Borwein
method proposed by Raydan in \cite{Raydan}.

\subsection{Newton's method}

The basic idea of Newton's method for unconstrained optimization is to iteratively
use the quadratic approximation to the objective function $f$ at the current point $\bold x_k$ and to minimize the approximation.
It is second order optimization method that is given by iterative rule

\begin{equation}\label{Newton_method}
\bold x_{k+1} = \bold x_k - t_k G_k^{-1}\bold g_k, \quad k = 0,1,\ldots.
\end{equation}

Clearly, Newton's method is called second order method, since it uses the information obtained from the second
partial derivatives of the objective $f$ which are incorporated in the Hessian $G(\bold x)$. Newton's direction
is descent direction if and only if the Hessian is positive definite matrix. It is much more faster than first order
methods, but it suffers from both expensive Hessian computation and inversion as well as constrain that the Hessian
need to be positive definite. Usually, it is not applicable to large scale problems.

\subsection{Conjugate Gradient methods}

The conjugate gradient method is an algorithm that was design as a numerical solution of particular systems of linear equations,
namely those whose matrix is symmetric and positive-definite. It was originally developed by Hestenes and Stiefel \cite{Hestenes}.
Later, the conjugate gradient method is adopted and modified and nowadays can be used to solve unconstrained optimization problems.
\smallskip

The first conjugate gradient method for solving unconstrained optimization problems was introduced by Fletcher and Reeves \cite{Fletcher}.
It is one of the earliest known techniques for solving large-scale nonlinear optimization problems.
Over the years, many variants of this original scheme have been introduced and some are widely used in practice.
The key features of these algorithms are that they require no matrix storage and are much faster than the original first
order methods such as gradient descent method. Conjugate gradient methods follows the iterative scheme
\begin{equation}\label{conjugate_gradient}
\bold x_{k+1} = \bold x_k + t_k \bold p_k, \quad k = 0,1,\ldots
\end{equation}
where new search direction $\bold p_{k+1}$ is determined as follows
\begin{equation}\label{cg_direction}
\bold p_{k+1} = -\bold g_{k+1} + \beta_{k+1} \bold p_k, \quad k = 0,1,\ldots
\end{equation}
The way to compute scalar value $\beta_k$ (which is essential for determining new conjugate direction $\bold p_k$)
is what distinguishes different modifications of the originally proposed method.
It is known (easy to check) that different formulae for computing scalar $\beta_k$ yield to the methods that are equivalent in the sense
that all produce the same search directions when used in minimizing a quadratic function with positive definite Hessian matrix.
However, for a general nonlinear function with inexact line search, their behavior is markedly different.
The strong Wolfe conditions are usually used in the analyses and implementations of various conjugate gradient methods.
Some descriptions will be given later in this subsection.
We present several the most famous and most widely used conjugate gradient methods. Conjugate gradient methods are known
as a methods with super linear rate of convergence.

\vskip 8pt
{\bf \noindent Fletcher-Reeves}
\vskip 5pt

Using the fact that the solving a linear system is equivalent to minimizing a positive definite quadratic
function, Fletcher and Reeves \cite{Fletcher} in the 1960s modified original conjugate gradient method and developed
a conjugate gradient method for unconstrained minimization. They introduce the following formula for the unknown
parameter
$$\beta_k = \frac{\bold g^T_k \bold g_k}{\bold g^T_{k-1} \bold g_{k-1}}.$$
Numerical properties showed the the given method outperform the well known steepest descent method.

\vskip 8pt
{\bf \noindent Polak-Ribiere}
\vskip 5pt

The next modification is introduced by Polak and Ribiere \cite{Polak}, which defines the parameter $\beta_k$ as follows

$$\beta_k = \frac{\bold g_k^T \bold y_{k-1}}{\left\lVert \bold g_{k-1} \right\rVert^2}.$$

Numerical experience indicates that this algorithm tends to be the more robust and efficient than Fletcher-Reeves method.
Surprisingly, strong Wolfe conditions do not guarantee that the direction $\bold p_k$ is always a descent direction.
With a simple adaptation of the parameter $\beta_k$

$$\beta_k^+ = \max( \beta_k, 0) $$
a strong Wolfe conditions ensures that the descent property holds.

\vskip 8pt
{\bf \noindent Hestenes-Stiefel}
\vskip 5pt

Another modification that coincide with the Fletcher-Reeves formula is proposed by Hestenes and Stiefel, see \cite{Hestenes}.
Their choice for parameter $\beta_k$ defined by
$$\beta_k = \frac{\bold g^T_k \bold y_{k-1}}{\bold p^T_{k-1}\bold y_{k-1}},$$
gives rise to an algorithm that is similar to Polak-Ribiere method, both in terms of its theoretical
convergence properties and in its practical performance.

\vskip 21pt
{\bf \noindent Dai-Yuan}
\vskip 5pt

Dai and Yuan, in their paper \cite{DaiYuan}, presented new formula for the parameter $\beta_k$, introducing another conjugate gradient method
$$\beta_k = \frac{\bold g^T_k \bold g_k}{\bold p^T_{k-1} \bold y_{k-1}}.$$
It is shown that this new method is globally convergent as long as the standard Wolfe conditions
(not necessarily strong Wolfe) are satisfied.
Moreover, they claimed that the conditions on the objective function are also weaker than the usual ones.

\vskip 8pt
{\bf \noindent CG\_Descent}
\vskip 5pt

A new nonlinear conjugate gradient method are proposed and analyzed by Hager and Zhang \cite{cg_descent1, cg_descent2}.
Their choice for parameter $\beta_k$ is defined as follows
$$\beta_k = \frac{1}{p_y} \left(\bold y_{k-1}-\frac{2 \|\bold y_{k-1}\|^2}{p_y} \bold p_{k-1}\right)^T \bold g_k, \quad
\rm{where}\ \ p_y = \bold p_{k-1}^T \bold y_{k-1}.$$
\smallskip

A global convergence result for CG\_Descent is established when the line search fulfills the Wolfe conditions. Additionally, new line search algorithm
is developed that is more efficient and highly accurate. High accuracy is achieved by using a
convergence criterion called {\em approximate Wolfe} conditions, obtained by replacing the
sufficient decrease criterion in the Wolfe conditions with an approximation that can be evaluated with
greater precision in a neighborhood of a local minimum. 
Presented numerical results for new method as well as new line search are given and the authors claim that their method
outperform other conjugate gradient methods on standard unconstrained optimization problems.

\subsection{Modified Newton's methods}

In order to use nice properties of original Newton's method and avoid some of it's drawbacks a new class of optimization
methods is introduced, so called modified Newton methods (also known, in literature, as Newton's like methods).
This group includes methods that are developed as a solution which
uses modified Newton's direction. Namely, following the appropriate ideas these methods combine Hessian and gradient in order
to dynamically compute the search direction which will provide the sufficient function decrease. Three different solvers that belongs
to this group are presented.

\vskip 8pt
{\bf \noindent Goldstein-Price}
\vskip 5pt

Goldstein and Price \cite{GoldsteinPrice} presented one modification of Newton's method, also studied in \cite{SY}.
The general idea was to use the steepest descent direction $\bold d_k = −\bold g_k$ in situations when $G_k$ is not positive definite,
and thus not guarantee the convergence.
They analysed so called angle rule
$$ \theta \leq \frac{\pi}{2} - \mu,\quad \rm{for\ some\ small}\quad \mu > 0,$$
where $\theta$ present the angle between the vectors $\bold d_k$ and $-\bold g_k$.
In order to satisfied the descent condition, the angle rule for search direction $\bold d_k$ should be satisfied.
Taking this into account they proposed the following choice for computing search direction
$$
\bold d_k =
\begin{cases}
-G_k^{-1} \bold g_k,& \text{if}\ \cos\theta\geq \eta,\\  -\bold g_k, &\text{otherwise.}
\end{cases}
$$
In our implementation $\eta = 0.2$ is used. Some other values for this parameter can be also considered.

\vskip 8pt
{\bf \noindent Levenberg-Marquardt}
\vskip 5pt

Levenberg-Marquardt method is originally constructed and is well studied as a solver for non-linear least squares problems.
Levenberg-Marquardt method is obtained as a generalization of Gauss-Newton methods by replacing the
line search strategy with trust region strategy. Trust region strategy is imposed to avoids one of the weaknesses of Gauss-Newton method,
its behavior when the Jacobian $J(\bold x)$ is rank deficient, or nearly so.
This method can be very easily modified in order to serve as a solver for unconstrained optimization non-linear problems.
Namely, instead of $J^TJ$ which represent the Hessian approximation (for least squares problems) the
true Hessian can be used in the case of minimizing general nonlinear objective functions.
\smallskip

By applying already mentioned slight change the search direction can be computed as follows
$$
\bold d_k = -(G_k + \lambda I)^{-1} \bold g_k, \quad \lambda > 0.
$$
This term for computing the search direction is well known as Levenberg method, originally proposed in \cite{Levenberg}.
Similar idea is used in \cite{Goldfeld} in order to determine search direction by using information obtained from Hessian.
More precisely, when Hessian is not positive definite,
one changes the model Hessian into $G_k + \lambda I$ such that $G_k + \lambda I$ becomes positive definite.
There is another slight modification of Levenberg search direction proposed by Marquardt
$$
\bold d_k = -\left[G_k + \lambda \cdot diag(G_k)\right]^{-1} \bold g_k, \quad \lambda > 0,
$$
also known as Levenberg-Marquardt method \cite{Marquardt}. It is claimed that this slight modification
produces some small numerical improvements. These methods can be viewed as a switch rule between Gauss-Newton method
and the steepest descent method. By changing the values of parameter $\lambda$ we can chose any direction between these
two directions. In the case $\lambda = 0$ it becomes Gauss-Newton direction, while putting $\lambda$ to be very large
it approximate steepest descent direction. In each iteration algorithm search for best value of $\lambda$ which
sufficiently decrease the function value.

\subsection{Quasi-Newton methods}

Quasi-Newton method is a class of methods which need not to compute the Hessian, but generate a series
of Hessian approximations, and at the same time maintain a fast rate of convergence.
This method use function value and its gradient during the iteration process as well as corresponding
Hessian (or it's inverse) approximation. By this idea, expensive Hessian and its inverse computation is avoided.
On the other side, by appropriate choice of Hessian approximation good convergence
rate is retained. For the sake of simplicity, we use the following notations for Hessian and Hessian inverse approximation,
$B_k$ and $H_k$. Methods from this group satisfy so called secant (quasi-Newton) equation
\begin{equation} \label{secantEq}
B_{k+1}\bold s_k = \bold y_k.
\end{equation}
Additionally, dual equation can be very easily derived
\begin{equation} \label{secantEqDual}
H_{k+1}\bold y_k = \bold s_k. \end{equation}
Four different methods are currently presented.

\vskip 8pt
{\bf \noindent SR1 approximation}
\vskip 5pt

A simple symmetric rank one update (SR1) that satisfies secant equation \eqref{secantEqDual} is given by the following equation

\begin{equation} \label{SR1}
H_{k+1}=H_k+\frac{(\bold s_k-H_k\bold y_k)(\bold s_k-H_k\bold y_k)^T}{(\bold s_k-H_k\bold y_k)^T \bold y_k}.
\end{equation}

This general rank one update invented by Broyden in \cite{Broyden1} was originally developed for solving systems of nonlinear equations.
One of the main drawback of this update formula is that does not retain the positive definiteness of $H_k$.
Another issue is possibility of the denominator $(\bold s_k-H_k\bold y_k)^T \bold y_k$ to be very small or zero, which results in serious
numerical instability. In order to overcome some of this difficulties the following simple modification
(see \cite{Nocedal}, \cite{SY}) is proposed
$$
H_{k+1} =
\begin{cases}
\text{SR1\ update},& \text{if}\quad |(\bold s_k-H_k\bold y_k)^T \bold y_k| \geq r \|\bold s_k-H_k\bold y_k\| \| \bold y_k\|,\\
H_k, &\text{otherwise.}
\end{cases}
$$
where $r$ is some small number $r\in(0, 1)$, say $r = 10^{-8}$. The SR1 update formula is usually combined with the trust region methods,
but can be also accompanied by the line search algorithms.

\vskip 8pt
{\bf \noindent DFP}
\vskip 5pt

This method follows the quasi-Newton update proposed originally by
Davidon \cite{Davidon} and developed later by Fletcher and Powell \cite{FletcherPowell};
thus it is called Davidon-Fletcher-Powell formula (or DFP method).

\begin{equation} \label{DFP update}
H_{k+1}^{DFP}=H_k+\frac{\bold s_k\bold s_k^T}{\bold s_k^T\bold y_k}-\frac{H_k\bold y_k\bold y_k^TH_k}{\bold y_k^TH_k\bold y_k}.
\end{equation}
The DFP formula  finds the solution to the secant equation \eqref{secantEqDual} that is closest to the current estimate and satisfies
the curvature condition. It was the first quasi-Newton method to generalize the secant method to a multidimensional
problem. This update maintains the symmetry and positive definiteness of the Hessian approximation, which is essential for
convergence properties. The DFP method was very popular rank two update quasi Newton method, quite effective, but it was soon
superseded by the so called BFGS formula, which is its dual update. The DFP method as well as BFGS method achieved very good
results after applying the Wolfe conditions on step-size computation.
\vskip 8pt
{\bf \noindent BFGS}
\vskip 5pt

Another very important algorithm discovered independently by Broyden \cite{Broyden},
Fletcher \cite{FletcherBFGS}, Goldfarb \cite{Goldfarb} and Shanno \cite{Shanno} is named
Broyden-Fletcher-Goldfarb-Shanno (or shortly BFGS) method.
The BFGS method is one of the most popular members of this class.
Instead of imposing conditions on the Hessian inverse approximations $H_k$ (like for DFP method),
similar conditions on Hessian approximation $B_k$ are impose which yields to the following formula
\begin{equation} \label{BFGS update}
B_{k+1}^{BFGS}=B_k+\frac{\bold y_k\bold y_k^T}{\bold y_k^T\bold s_k}-\frac{B_k\bold s_k\bold s_k^TB_k}{\bold s_k^TB_k\bold s_k},
\end{equation}
which satisfies secant equation given by \eqref{secantEq}.
In order to obtain more efficient formula for updating $H_k$ instead of $B_k$, Sherman-Morrison-Woodbury formula
is applied. One of the benefit of this method is the fact that the BFGS formula (unlike DFP) has very
effective self-correcting properties. This is very important property in situations in which the BFGS can
produce bad results, such as incorrect estimates of the curvature of the objective function. The self correcting
properties of BFGS hold only when an adequate line search (like the Wolfe line search) is performed.

\vskip 8pt
{\bf \noindent L-BFGS}
\vskip 5pt
One of the main drawback of previously described quasi-Newton methods is it's inefficiency
in solving large scale problems. Namely, Hessian matrix approximation cannot be computed at a
reasonable cost, it is too expensive to store the whole matrix, and it's not easy to manipulate with.
In order to overcome these difficulties
a number of methods has been introduced. One of the most famous is so called Limited-memory BFGS
(shortly L-BFGS) proposed by Liu and Nocedal \cite{l-bfgs}. This method, which its name suggests, is based on
the BFGS updating formula. Instead of storing fully dense $n \times  n$ matrix, only a few vectors of length
$n$ that represent the approximations are stored. Despite these optimal storage requirements,
it yields good rate of convergence. In order to determine Hessian approximation, the curvature information
are used from only the few recent iterations. Curvature information which can be obtained from earlier
iterations is discarded.
\smallskip

The L-BFGS shares many features with other quasi-Newton methods, but is very different in how the
matrix-vector multiplication for finding the search direction $\bold d_k = H_k \bold g_k$ is carried out.
There are several approaches using a history of updates to form this direction vector.
One of the most common approach is the so-called two loop recursion, see \cite{Nocedal}.
In order to guarantee that the curvature condition will be satisfied the Wolfe line search is suggested to be used.
As far as we know the most optimal choice for the line search is More-Thuente line search algorithm (see \cite{More-Thuente})
which is the default choice in our application.

\subsection{Trust region methods}

In the trust region strategy, the main idea is to use the objective function $f(x)$ to construct it's approximation $m_k(x)$
whose behavior near the current point $\bold x_k$ is similar to that of the actual function $f(x)$. The function $m_k$ represent the
model function in $k$-th iteration. Instead of looking for the minimizer of the objective $f$ the minimizer of model function
$m_k$ (usually more simplify function than original one) will be determined. To avoid the situation in which the model $m_k$
may not be a good approximation of $f$ (optimal point $\bold x$ is far from current $\bold x_k$) the search for a minimizer
is restricted to some region around $x_k$. The region around the point $x_k$ is known as the trust region.
The model function is usually given as a quadratic model
\begin{equation} \label{trustModel}
m_k(\bold x_k + \bold d) = f_k + \bold d^T \bold g_k + \frac{1}{2}\bold d^T B_k \bold d,
\end{equation}
where $B_k$ is either Hessian or its approximation. The idea is to find vector $\bold d$  by
solving the following subproblem
\begin{equation} \label{trustProblem}
\min_{\bold d} m_k(\bold d) = \min_{\bold d} f_k + \bold d^T \bold g_k + \frac{1}{2}\bold d^T B_k \bold d,
\quad {\rm s.t.} \ \|\bold d\| \leq \Delta,
\end{equation}
where $\Delta$ is given trust region.

\vskip 8pt
{\bf \noindent Dogleg}
\vskip 5pt

One of the first idea for solving the trust-region subproblem \eqref{trustProblem} is the
so called dogleg method introduced by Powell \cite{Powell_Dogleg}.
To find an approximate solution of the subproblem \eqref{trustProblem}, inside the trust region,
Powell used a path consisting of two line segments. The first line segment runs from the current point $\bold x_k$
to the Cauchy point (a minimizer generated by the steepest descent method). Additionally, the second
line segment runs from the Cauchy point to the Newton point (the minimizer  generated by Newton method).
Let we denote the Cauchy point $P^C$ and Newton point $P^B$, then the path formally can be
defined as follows
\begin{equation}
p(\tau) =
\begin{cases}
\tau P^C,& 0 \leq \tau \leq 1,\\
P^C + (1-\tau)(P^B - P^C), &1 \leq \tau \leq 2.
\end{cases}
\end{equation}
The dogleg method determines $\tau$ in order to minimize the model $m_k$ along this path, inside the trust region.
In fact, it is not even necessary to carry out a search, because the
dogleg path intersects the trust region boundary at most once. The model function $m_k$ decreases along the path and it is
proved that the intersection point can be computed analytically.

\vskip 8pt
{\bf \noindent Dogleg-SR1}
\vskip 5pt

Instead of using the original Newton point $P^B = \bold x_k - G_k^{-1}\bold g_k$ in the dogleg method, it is possible
to use point obtained after applying some of the Newton inverse approximation. Namely, some ideas from quasi-Newton
methods can be established. The most popular idea is to use SR1 quasi-Newton approximation based on its ability to generate very
good Hessian approximations. Thus, the Newton point is obtained as follows
\begin{equation}
P^B = \bold x_k - H_k \bold g_k,
\end{equation}
where matrix $H_k$ is determined by  \eqref{SR1}.
Additionally, the subproblem given by \eqref{trustProblem} is determined by SR1 quasi Newton
approximation $H_k^{-1}$. This method is pretty suitable in the case of difficult and expensive Hessian calculation.
To obtain a fast rate of convergence, it is important for the matrix $H_k$ to be updated
in every iteration. Namely, along the failed direction in which the step was poor, $H_k$ should be updated because it represents
the poor approximation of the true Hessian in this direction. The idea for the implementation of this method
is taken from \cite{Nocedal}.

\subsection{Line search methods}\label{lineSearchList}

In this subsection we introduce the list of available line searches as well as some very shorts explanations.
Detailed information about the specific line search can be found in original papers as well as in some relevant
book or monograph, such as  \cite{FletcherBook}, \cite{Nocedal} and \cite{SY}.
\smallskip

Just to emphasize that the most of the parameters that corresponds to the specific line search algorithm can by 
manually tuned throughout the application interface under the panel '\texttt{Line search params}'.

\vskip 8pt
{\bf \noindent Fixed Step Size}
\vskip 5pt

Fixed step size is the simplest line search rule. It simply reads value for parameter $t$ from application
interface. During the iterative process step size $t$ is fixed.
This line search is predefined as default choice for Newton line search method and is recommended by the authors
for this method.

\vskip 8pt
{\bf \noindent Step size determined by previous values }
\vskip 5pt

Next two methods follow the heuristics which use information from previous iterations. Step-size is computed dynamically
according to the function and step-size values taken from previous iterations.
First one named {\em correction by previous iteration} decreases step-size, by some factor $0< c_1 <1$. Namely, if current point
try to moves away from the local minimum, with respect to the previous point, step-size decreases, otherwise remains the same,
see algorithm \ref{twoHeuritics} (left).
Therefore, generally speaking, as we are closer to the solution, the step-size value reduces.
In our implementation $c_1 = 0.5$ is established; some other values can also be considered.
One of the drawback of this heuristic is the fact that the step-size can only be reduced. Thus, if the current point
is far from the solution, and initial step-size is not appropriate (much smaller than optimal) the convergence can be very slow.
\smallskip

Second one, named {\em correction by previous two iteration} is more robust and tries to resolve the mentioned issue of previous
heuristic. Namely, the current step-size is computed according to the function values in previous two iterations,
see algorithm \ref{twoHeuritics} (right).

\begin{multicols}{2}
[
\begin{algo}\label{twoHeuritics}
Step-size determined by previous iterations values
\end{algo}
]
\begin{small}
\begin{verbatim}
----------------------
CorrPrevIter
----------------------
while (f_{k+1} >= f_k)
    t = t * c1;
    k = k + 1;
    x_{k+1} = x_k + t*d_k
    determine f_{k+1}
end


----------------------
CorrPrevTwoIter
----------------------
if (f_{k+1} < f_k && f_k < f_{k-1})
    t = t * c2; % increase step length
else if f_{k+1} >= f_k
        % decrease step length
        t = t * c1;
    end
end
\end{verbatim}
\end{small}
\end{multicols}
\noindent If the current point is far from the solution and if in previous two iterations the function value decreases than the
step-size can be increased by some factor $c_2 > 1$. Otherwise, the same idea is applied as in previous line search heuristic.
In current implementation the following value $c_2 = 1.2$ for increasing factor is chosen.
From our point of view this two values best fit parameters $c_1$ and $c_2$ but also some other values can be cosidered.

\vskip 8pt
{\bf \noindent Backtracking}
\vskip 5pt

We continue with the one of the easiest inexact line search method also known as Armijo's (Backtracking) line search \cite{Armijo}.
This method is described by the following problem:
determine parameter $t > 0$ such that the following condition
\begin{equation}\label{Armijo}
f(\bold x_k+t \bold d_k)\leq f(\bold x_k)+\rho t \bold g_k^T \bold d_k,
\end{equation}
be satisfied, where $\rho \in (0,1)$.
We present two different implementations that satisfy Armijo condition \eqref{Armijo}, named Backtracking and Armijo.
Because of its simplicity, historical importance and practical value we decide to keep the implementation of the simpler idea.
The method known as Backtracking is still in use, see for example \cite{Andrei3}, \cite{Stanimirovic}.
The idea is to take initial value for $t$ (usually $t=1$) and
iterate (using backtracking algorithm) in order to make condition \eqref{Armijo} satisfied.
It finds $t$ by applying $t = t\cdot\beta$ until \eqref{Armijo} becomes satisfied, where $0<\beta<1$.
Therefore, the smallest integer $l$ needs to be determined such that following inequality
\begin{equation} \label{Backtracking}
f(\bold x_k + \beta^{l}t\bold d_k) \leq f(\bold x_k) + \rho\beta^{l}t\bold g_k^{T}\bold d_k, \quad l = 0,1,\ldots
\end{equation}
holds.

\vskip 8pt
{\bf \noindent Armijo line search}
\vskip 5pt

Another version, named Armijo, uses function interpolation to find the best possible step-size $t$ which satisfies \eqref{Armijo}.
The solution that uses interpolation gives better numerical performances with respect to the solution obtained by the simple backtracking.
For detailed explanation see, for example, \cite{Nocedal}.
\smallskip

Note that, in context of Armijo line search, parameters in application interface with names '\texttt{start point}', '\texttt{beta}' and '\texttt{rho}'
(under '\texttt{Line search params}' panel) correspond to $t$, $\beta$ and $\rho$ in \eqref{Armijo} and \eqref{Backtracking}.

\vskip 8pt
{\bf \noindent Goldstein line search}
\vskip 5pt

Next line search rule is Goldstein rule, established by its author in \cite{Goldstein}. This line search represent the
improvement of the previous Armijo rule. Namely, besides the upper bound for the step-size $t$ given by
condition \eqref{Armijo}, lower bound, given by \eqref{Goldstein2}, is introduced. Therefore, Goldstein rule is determined
by following two expressions
\begin{equation}\label{Goldstein1}
f(\bold x_k+t \bold d_k)\leq f(\bold x_k)+\rho t \bold g_k^T \bold d_k
\end{equation} and
\begin{equation}\label{Goldstein2}
f(\bold x_k+t \bold d_k)\geq f(\bold x_k)+(1-\rho)t \bold g_k^T \bold d_k,
\end{equation}
where $0< \rho <\frac 12$. Step-size $t$ that satisfies this two conditions is computed by binary search on starting interval.
The parameters '\texttt{start point}' and '\texttt{rho}' in application interface corresponds to $t$ and $\rho$ in \eqref{Goldstein1} 
and \eqref{Goldstein2} in the case when Goldstein line search is selected.

\vskip 8pt
{\bf \noindent Wolfe line search}
\vskip 5pt

In order to improve Goldstein rule another idea is introduced, named after its author, Wolfe line search.
Namely, there is no guarantee that the lower bound given
by \eqref{Goldstein2} will contain a local minimum. Therefore, Wolfe in his paper \cite{Wolfe} introduced additional
condition (given by \eqref{Wolfe}) as a lower bound for determining parameter $t$ for which has been proven to contain a local minimum.
Thus, step-size $t$ is chosen such that it satisfies Armijo rule \eqref{Goldstein1} and the following curvature condition
\begin{equation}\label{Wolfe}
\bold g_{k+1}^T \bold d_k \geq \sigma \bold g_k^T \bold d_k,\quad \sigma \in (\rho,1).
\end{equation}
Similarly as for previous line searches, parameter '\texttt{sigma}' from application interface corresponds
to $\sigma$ from \eqref{Wolfe} and \eqref{StrongWolfe}.

\vskip 8pt
{\bf \noindent Strong Wolfe line search}
\vskip 5pt

Additionally, Wolfe presented another slight modification of his original Wolfe line search, known as strong Wolfe.
In fact, instead of using condition \eqref{Wolfe} another stronger condition \eqref{StrongWolfe} is accompanied to the
Armijo rule \eqref{Goldstein1}. This condition ensures that interval in which the step-size need to be determined is more closer
to the one dimensional minimizer.
\begin{equation} \label{StrongWolfe}
|\bold g_{k+1}^T \bold d_k| \leq \sigma |\bold g_k^T \bold d_k|.
\end{equation}
\smallskip

For the implementation of Wolfe and strong Wolfe algorithms we chose so called two stage algorithm, see for example \cite{Nocedal}.
First stage begins with a trial estimate, and increases it until it finds an acceptable step-size or appropriate interval.
In the second stage function called zoom is used, which successively decreases the size of the interval until an
acceptable step length is identified.


\vskip 8pt
{\bf \noindent Approximate Wolfe line search}
\vskip 5pt

In order to get better numerical stability, the authors in \cite{cg_descent1, cg_descent2} proposed a new line search scheme
which is efficient and highly accurate.  Efficiency is achieved by exploiting properties of
linear interpolants in a neighborhood of a local minimizer.  High accuracy is achieved by using a
adopted Wolfe criterion \eqref{ApproxWolfe}, which they call the 'approximate Wolfe' conditions
\begin{equation}\label{ApproxWolfe}
(2 \rho - 1) \bold g_k^T\bold d_k \geq  \bold g_{k+1}^T \bold d_k \geq \sigma \bold g_k^T \bold d_k,
\quad \rho < \min(0.5, \sigma).
\end{equation}
This criterion is obtained by replacing the sufficient decrease condition \eqref{Goldstein2} with an approximation
that can be evaluated with greater precision. Namely, the precision is much higher in a neighborhood of a local minimum
than the usual precision obtained by sufficient decrease criterion given by Armijo rule.
This line search is default option for $CG\_Descent$ conjugate gradient method developed by the same authors.

\vskip 8pt
{\bf \noindent More-Thuente line search}
\vskip 5pt

More-Thuente line search is a line search procedure for computing
step-size parameter such that it satisfies so called strong Wolfe condition \eqref{StrongWolfe}.
It is an iterative method which is proven to be very effective. Namely, the proposed method start with initial step-size $t_0 \in [\alpha_{min}, \alpha_{max}]$
and iteratively generates the sequence of nested intervals $\{I_k\}$ and sequence of possible step lengths $t_k \in I_k \cap [\alpha_{min}, \alpha_{max}]$
until it finds the one that satisfies strong Wolfe condition. The authors claim that the algorithm
terminates within a small number of iterations. It is one of the most popular method in line search category. This method is
originally proposed by J.J. More and D.J. Thuente, see \cite{More-Thuente}. This line search is set as a default line search for
L-BFGS method which is one of the state of the art algorithms for solving large scale unconstrained problems.

\vskip 8pt
{\bf \noindent Non-monotone line search}
\vskip 5pt

In NonMonotone line search strategy, introduced by Grippo et al in \cite{Grippo}, the condition that the function value decreases in each iteration
is not imposed. The Nonmonotone line search is based on the usage of a positive constant integer $M$.
In each iteration the step-size $t$ is obtained in such a manner to fulfil the inequality
\begin{equation} \label{non-monotone}
f(\bold x_k + t \bold d_k) \leq \max\limits_{0\leq j \leq m(k)} \left[ f(\bold x_{k-j})+\rho t \bold g_k^T \bold d_k\right],
\end{equation}
where $m(0)=0$,  $0 \leq m(k) \leq \min \{m(k-1)+1,M-1\}$, and $\rho$ is a parameter from the Armijo's rule \eqref{Armijo}.
\smallskip

It is clear that this line search can bee seen as a generalization of Armijo's line search.
The parameter '\texttt{M}' in application interface corresponds to $M$ in \eqref{non-monotone} in the case
when NonMonotone line search is selected. This line search is default option for Barzilai-Borwein method
as well as Scalar-Correction method as it is suggested by the authors of the proposed methods, see \cite{Miladinovic} and \cite{Raydan}.

\end{appendices}

\begin{footnotesize}

\end{footnotesize}

\end{document}